\documentclass[a4paper,11pt]{article}
\usepackage{amsfonts}
\usepackage{pstricks}
\usepackage{pst-node,pst-tree,pst-eps,xspace}
\IfFileExists{pst-beta.tex}
             {\input{pst-beta.tex}}
             {\RequirePackage{pst-tree}}

\usepackage{latexsym}
\usepackage{amsthm}
\usepackage{fancyheadings}
\usepackage{amsmath}
\usepackage{graphicx}
\usepackage{float}
\usepackage{authblk}
\usepackage[font=small]{caption}
\restylefloat{table}

\def\preccdot{\prec\hspace{-1.9mm}\cdot}

\newtheorem{Co}{Corollary}

\newtheorem{The}{Theorem}
\newtheorem{Pro}{Proposition}
\theoremstyle{definition}
\newtheorem{De}{Definition}

\newcommand{\bsb}[1]{\boldsymbol{#1}}

\title{More restricted growth functions: Gray codes and exhaustive generations}
\author[1]{Ahmad Sabri} 
\author[2]{Vincent Vajnovszki}
\affil[1]{Center for Computational Mathematics Studies, Gunadarma University\authorcr \texttt{sabri@staff.gunadarma.ac.id}}
\affil[2]{LE2I, Universit\'e Bourgogne Franche-Comt\'e\authorcr \texttt{vvajnov@u-bourgogne.fr}}
\begin{document}
\maketitle
\begin{abstract}
A Gray code for a combinatorial class is a method for listing the objects in the class 
so that successive objects differ in some prespecified, small way, typically expressed 
as a bounded Hamming distance.
In a previous work,
the authors of the present paper showed, among other things, 
that the $m$-ary Reflected Gray Code Order yields a Gray code for the set of restricted growth 
functions.
Here we further investigate variations of this order relation,
and give the first Gray codes and efficient generating algorithms
for bounded restricted growth functions.


\end{abstract}

\section{Introduction}

In \cite{SV} the authors shown that both the order relation induced by the generalization
of the Binary Reflected Gray Code and one of its suffix partitioned version yield Gray codes 
on some sets of restricted integer sequences, and in particular for restricted growth functions.
These results are presented in a general framework,
where the restrictions are defined by means of statistics on integer sequences.

In the present paper we investigate two prefix partitioning order relations on the set 
of {\it bounded} restricted growth functions: as in \cite{SV}, the original Reflected Gray Code Order
on $m$-ary sequences, and a new order relation which is an appropriate modification of the 
former one. We show that, according to the parity of the imposed bound, one of these order relations
gives a Gray code on the set of bounded restricted growth functions. 
As a byproduct, we obtain a Gray code for restricted growth functions with a specified odd value for the largest entry;
the case of an even value of the largest entry remains an open problem.
In the final part we present the corresponding 
exhaustive generating algorithms.
A preliminary version of these results were presented at
The Japanese Conference on Combinatorics and its Applications in May 2016 in Kyoto \cite{SV_Jap}.

\section{Notation and definitions}

A {\it restricted growth function} of length $n$ is an integer sequence 
$\bsb{s}=s_1s_2\ldots s_n$ with $s_1=0$ and $0\leq s_{i+1}\leq \max\{s_j\}_{j=1}^i+1$, for all $i$, $1\leq i\leq n-1$. 
We denote by $R_n$ the set of length $n$ restricted growth functions,
and its cardinality is given by the $n$th Bell number (sequence {\rm A000110} in \cite{sloa}), 
with the exponential generating function $e^{e^x}-1$.
And length $n$ restricted growth functions encode the partitions of an $n$-set.

For an integer $b\geq 1$, let $R_n(b)$ denote the set of {\it $b$-bounded} sequences in $R_n$, that is, 
$$R_n(b)=\{s_1s_2\ldots s_n\in R_n\,:\, \max\{s_i\}_{i=1}^n\leq b\},
$$
and
$$R_n^{*}(b)=\{s_1s_2\ldots s_n\in R_n\,:\, \max\{s_i\}_{i=1}^n= b\}.
$$
See Table \ref{Tb1} for an example.

\begin{table}[h]
\centering
\caption{\label{Tb1}
The set $R_5(2)$, and in bold-face the set $R^*_5(2)$.
Sequences are listed in $\preccdot$ order (see Definition \ref{de:coRGCorder}) and in italic is the Hamming distance 
between consecutive sequences.}
\begin{tabular}{|rcc|rcc|rcc|}
\hline      
1.  &  0\, 0\, 0\, 0\, 0  &         & 15. & 0\, 1\, 0\, 0\, 0 & {\it 3} & 29. & \bf{0 1 1 1 2} & {\it 1}\\
2.  &  0\, 0\, 0\, 0\, 1  & {\it 1} & 16. & 0\, 1\, 0\, 0\, 1 & {\it 1} & 30. & \bf{0 1 1 2 2} & {\it 1}\\
3.  &  0\, 0\, 0\, 1\, 0  & {\it 2} & 17. & \bf{0 1 0 0 2}    & {\it 1} & 31. & \bf{0 1 1 2 1} & {\it 1}\\
4.  &  0\, 0\, 0\, 1\, 1  & {\it 1} & 18. & 0\, 1\, 0\, 1\, 0 & {\it 2} & 32. & \bf{0 1 1 2 0} & {\it 1}\\
5.  &\bf{0 0 0 1 2}       & {\it 1} & 19. & 0\, 1\, 0\, 1\, 1 & {\it 1} & 33. & \bf{0 1 2 2 0} & {\it 1}\\
6.  &  0\, 0\, 1\, 0\, 0  & {\it 3} & 20. & \bf{0 1 0 1 2}    & {\it 1} & 34. & \bf{0 1 2 2 1} & {\it 1}\\
7.  &  0\, 0\, 1\, 0\, 1  & {\it 1} & 21. & \bf{0 1 0 2 2}    & {\it 1} & 35. & \bf{0 1 2 2 2} & {\it 1}\\
8.  &\bf{0 0 1 0 2}       & {\it 1} & 22. & \bf{0 1 0 2 1}    & {\it 1} & 36. & \bf{0 1 2 1 2} & {\it 1}\\
9.  &  0\, 0\, 1\, \,1 0  & {\it 2} & 23. & \bf{0 1 0 2 0}    & {\it 1} & 37. & \bf{0 1 2 1 1} & {\it 1}\\
10. &  0\, 0\, 1\, 1\, 1  & {\it 1} & 24. & 0\, 1\, 1\, 0\, 0 & {\it 2} & 38. & \bf{0 1 2 1 0} & {\it 1}\\
11. &\bf{0 0 1 1 2}       & {\it 1} & 25. & 0\, 1\, 1\, 0\, 1 & {\it 1} & 39. & \bf{0 1 2 0 2} & {\it 2}\\
12. &\bf{0 0 1 2 2}       & {\it 1} & 26. & \bf{0 1 1 0 2}    & {\it 1} & 40. & \bf{0 1 2 0 1} & {\it 1}\\
13. &\bf{0 0 1 2 1}       & {\it 1} & 27. & 0\, 1\, 1\, 1\, 0 & {\it 2} & 41. & \bf{0 1 2 0 0} & {\it 1}\\
14. &\bf{0 0 1 2 0}       & {\it 1} & 28. & 0\, 1\, 1\, 1\, 1 & {\it 1} & & &\\
\hline  
\end{tabular}
\end{table}

If a list of same length sequences is such that the Hamming distance between successive sequences 
(that is, the number of positions in which the sequences differ) is bounded from above by a constant, 
independent on the sequences length, then the list is said to be a {\it Gray code}. 
When we want to explicitly specify this constant, say $d$, then we refer to such a list as a {\it $d$-Gray code};
in addition, if the positions where the successive sequences differ are adjacent,
then we say that the list is a {\it $d$-adjacent Gray code}.

The next two definitions give order relations on the set of $m$-ary integer sequences of 
length $n$ on which our Gray codes are based.

\begin{De}
\label{de:RGCorder}
Let $m$ and $n$ be positive integers with $m\geq 2$. 
The {\em Reflected Gray Code Order} $\prec$ on $\{0,1,\ldots,m-1\}^n$ is defined as:
$\bsb{s}=s_1s_2\ldots s_n$ is less than $\bsb{t}=t_1t_2\ldots t_n$,
denoted by $\bsb{s}\prec\bsb{t}$, if

\begin{center}
 either $\sum_{i=1}^{k-1} s_i$ is even and $s_k<t_k$, or $\sum_{i=1}^{k-1} s_i$ is odd and $s_k>t_k$
\end{center}
for some $k$ with $s_i=t_i$ ($1\leq i\leq k-1$) and $s_k\neq t_k$.
\end{De}

This order relation is the natural extension to $m$-ary sequences of the order 
induced by the Binary Reflected Gray Code introduced in \cite{Gray}.
See for example \cite{BBPSV,SV} where this order relation and its variations
are considered in the context of factor avoiding words and of statistic-restricted sequences.

\begin{De}
\label{de:coRGCorder}
Let $m$ and $n$ be positive integers with $m\geq 2$.
The {\em co-Reflected Gray Code Order}\footnote
{In \cite{SV} a similar terminology is used for a slightly different notion}
$\preccdot$ on $\{0,1,\ldots,m-1\}^n$ is defined as:
$\bsb{s}=s_1s_2\ldots s_n$ is less than $\bsb{t}=t_1t_2\ldots t_n$,
denoted by $\bsb{s}\preccdot\ \bsb{t}$, if 
\begin{center}
either  $U_k$ is even and $s_k<t_k$, or $U_k$ is odd and $s_k>t_k$
\end{center}
for some $k$ with $s_i=t_i$ ($1\leq i\leq k-1$) and $s_k\neq t_k$ where 
$
U_k= | \{i\in \{1,2,\ldots,k-1\}:s_i\neq 0, s_i \mbox{ is even}\}|$.
\end{De}
See Table \ref{Tb_co} for an example.

For a set $S$ of same length integer sequences the {\it $\prec$-first} (resp. {\it $\prec$-last}) sequence 
in $S$ is the first (resp. last) sequence when the set is listed  in $\prec$ order; and {\it $\preccdot$-first} 
and {\it $\preccdot$-last} are defined in a similar way. And for sequence $\bsb u$, $ \bsb u\,|\,S$ denotes 
the subset of $S$ of sequences having prefix~$\bsb u$.

Both order relations, Reflected and co-Reflected Gray Code Order produce prefix partitioned lists,
that is to say, if a set of sequences is listed in one of these order relations, then the sequences having a common prefix are consecutive in the list.

\begin{table}[h]
\centering
\caption{\label{Tb_co}
The set $\{0,1,2\}^3$ listed in 
$\preccdot$ order.}
\begin{tabular}{|cc|cc|cc|}
\hline      
1.  & 0 0 0 & 10. & 1 0 0 & 19. & 2 2 0 \\ 
2.  & 0 0 1 & 11. & 1 0 1 & 20. & 2 2 1 \\ 
3.  & 0 0 2 & 12. & 1 0 2 & 21. & 2 2 2 \\  
4.  & 0 1 0 & 13. & 1 1 0 & 22. & 2 1 2 \\  
5.  & 0 1 1 & 14. & 1 1 1 & 23. & 2 1 1 \\ 
6.  & 0 1 2 & 15. & 1 1 2 & 24. & 2 1 0 \\  
7.  & 0 2 2 & 16. & 1 2 2 & 25. & 2 0 2 \\  
8.  & 0 2 1 & 17. & 1 2 1 & 26. & 2 0 1 \\  
9.  & 0 2 0 & 18. & 1 2 0 & 27. & 2 0 0 \\  
\hline  
\end{tabular}
\end{table}

\section{The Gray codes}

In this section we show that the set $R_n(b)$, with $b$ odd, listed in $\prec$ order is a Gray code.  
However, $\prec$ does not induce a Gray code when $b$ is even: 
the Hamming distance between two consecutive sequences can be arbitrary large for 
large enough $n$. 
To overcome this, we consider $\preccdot$ order instead of $\prec$ order
when $b$ is even, and we show that the obtained list is a Gray code.

In the proof of Theorem \ref{k_odd} below we need the following propositions which give
the forms of the last and first sequence in $R_n(b)$ having a certain fixed prefix, when sequences are listed in $\prec$
order.

\begin{Pro}
\label{pro:pro_Rb_odd1}
Let $b\geq1$ and odd, $k\leq n-2$ and $\bsb s=s_1\ldots s_{k}$. If $\bsb t$ is the $\prec$-last 
sequence in $\bsb s\,|\,R_n(b)$, then $\bsb t$ has one of the following forms:
\begin{enumerate}
\item $\bsb t=\bsb sM0\ldots0$ if $\sum_{i=1}^{k} s_i$ is even and $M$ is odd,
\item $\bsb t=\bsb sM(M+1)0\ldots0$ if $\sum_{i=1}^{k} s_i$ is even and $M$ is even,
\item $\bsb t=\bsb s0\ldots0$ if $\sum_{i=1}^{k} s_i$ is odd,
\end{enumerate}
where $M=\min\{b,\max\{s_i\}_{i=1}^k+1\}$.
%where $M=\max\{a_i\}_{i=1}^k+1$.
\end{Pro}
\proof
Let $\bsb t=s_1\dots s_kt_{k+1}\ldots t_n$ be the $\prec$-last sequence in $\bsb s\,|\,R_n(b)$.

\noindent
Referring to the definition of $\prec$ order in Definition \ref{de:RGCorder}, if 
$\sum_{i=1}^{k} s_i$ is even, then $t_{k+1}=\min\{b,\max\{s_i\}_{i=1}^k+1\}= M$, and 
based on the parity of $M$, two cases can occur.
\begin{itemize}
\item If $M$ is odd, then we have that $\sum_{i=1}^k s_i+t_{k+1}=\sum_{i=1}^{k} s_i+M$ is odd, thus 
      $t_{k+2}\ldots t_n=0\ldots 0$, and we retrieve the form prescribed by the first point of the proposition.
\item If $M$ is even, then $M\neq b$ and $\sum_{i=1}^k s_i+t_{k+1}=\sum_{i=1}^{k} s_i+M$ is even, thus 
      $t_{k+2}=\max\{s_1,\ldots, s_k,t_{k+1}\}+1=M+1$, which is odd. 
      Next, we have $\sum_{i=1}^k s_i+t_{k+1}+t_{k+2}=\sum_{i=1}^k s_i+2M+1$ is odd, and this implies as above that 
      $t_{k+3}\ldots t_n=0\ldots 0$, and we retrieve the second point of the proposition.
\end{itemize}

\noindent
For the case when $\sum_{i=1}^{k} s_i$ is odd, in a similar way we have $t_{k+1}\dots t_n=0\dots 0$.
\qed

\medskip
The next proposition is the `first' counterpart of the previous one.
Its proof is similar by exchanging the parity of the summation from `odd' to `even' and vice-versa,
and it is left to the reader.

\begin{Pro}
\label{pro:pro_Rb_odd2}
Let $b\geq1$ and odd, $k\leq n-2$ and $\bsb s=s_1\ldots s_{k}$. If $\bsb t$ is the $\prec$-first sequence in $\bsb s\,|\,R_n(b)$, 
then $\bsb t$ has one of the following forms:
\begin{enumerate}
\item $\bsb t=\bsb sM0\ldots0$ if $\sum_{i=1}^{k} s_i$ is odd and $M$ is odd, 
\item $\bsb t=\bsb sM(M+1)0\ldots0$ if $\sum_{i=1}^{k} s_i$ is odd and $M$ is even,
\item $\bsb t=\bsb s0\ldots0$ if $\sum_{i=1}^{k} s_i$ is even,
\end{enumerate}
where $M=\min\{b,\max\{s_i\}_{i=1}^k+1\}$.
%where $M=\max\{a_i\}_{i=1}^k+1$.
\end{Pro}

Based on Propositions \ref{pro:pro_Rb_odd1} and \ref{pro:pro_Rb_odd2}, we have the following theorem.

\begin{The} 
\label{k_odd}
For any $n,b\geq 1$ and $b$ odd, $R_n(b)$ listed in $\prec$ order
is a $3$-adjacent Gray code.
\end{The}
\proof
Let $\bsb s=s_1s_2\dots s_n$ and $\bsb t=t_1t_2\dots t_n$ be two consecutive sequences in $\prec$ ordered list for the set $R_n(b)$, 
with $\bsb s\prec \bsb t$, and let $k$ be the leftmost position where $\bsb s$ and $\bsb t$ differ. 
If $k\geq n-2$, then obviously $\bsb s$ and $\bsb t$ differ in at most three positions, otherwise let 
$\bsb s'=s_1\ldots s_k$ and $\bsb t'=t_1\ldots t_k$. 
Thus, $\bsb s$ is the $\prec$-last sequence in $\bsb s'\,|\,R_n(b)$ and $\bsb t$ is the $\prec$-first sequence in $\bsb t'\,|\,R_n(b)$. 
Combining Propositions \ref{pro:pro_Rb_odd1} and \ref{pro:pro_Rb_odd2} we have that, when $k\leq n-3$,
$s_{k+3}s_{k+4}\dots s_n=t_{k+3}t_{k+4}\dots t_n=00\dots 0$.
And since $s_i=t_i$ for $i=1\ldots,k-1$, the statement holds.
\qed

\medskip
Theorem \ref{k} below shows the Graycodeness of $R_n(b)$, $b\geq 1$ and even, listed in $\preccdot$ order, and as for 
Theorem \ref{k_odd} we need the next two propositions; in its proof we will make use of the Iverson bracket notation: 
$[P]$ is $1$ if the
statement $P$ is true, and 0 otherwise. 
Thus, for a sequence $s_1s_2\dots s_n$ and a $k\leq n$, 
$| \{i\in \{1,2,\ldots,k\}:s_i\neq 0 \mbox{ and } s_i \mbox{ is even}\}|=\sum_{i=1}^k[s_i\neq0 {\rm\ and\ } s_i {\rm\ is\ even}]$.

\begin{Pro}
\label{pro:pro_Rb_even1}
Let $b\geq2$ and even, $k\leq n-2$ and $\bsb s=s_1s_2\ldots s_k$. If $\bsb t$ is the $\preccdot$-last sequence in 
$\bsb s\,|\,R_n(b)$, 
then $\bsb t$ has one of the following forms:
\begin{enumerate}
\item $\bsb t=\bsb sM0\ldots0$ if $U_{k+1}$ is even and  $M$ is even,
\item $\bsb t=\bsb sM(M+1)0\ldots0$ if $U_{k+1}$ is even and $M$ is odd,
\item $\bsb t=\bsb s0\ldots0$ if $U_{k+1}$ is odd,
\end{enumerate}
where $M=\min\{b,\max\{s_i\}_{i=1}^k+1\}$ and $U_{k+1}=\sum_{i=1}^k[s_i\neq0 {\rm\ and\ } s_i {\rm\ is\ even}]$.
%where $M=\max\{s_i\}_{i=1}^k+1$.
\end{Pro}
\proof
Let $\bsb t=s_1\dots s_kt_{k+1}\dots t_n$ be  the $\preccdot$-last sequence in $\bsb s\,|\,R_n(b)$. 

\noindent
Referring to the definition of $\preccdot$ order in Definition \ref{de:coRGCorder}, if $U_{k+1}$ is even, then $t_{k+1}=\min\{b,\max\{s_i\}_{i=1}^k+1\}=M>0$, and based on the parity of $M$, two cases can occur.
\begin{itemize}
\item If $M$ is even, then $U_{k+1}+[t_{k+1}\neq0 {\rm\ and\ } t_{k+1} {\rm\ is\ even}]=U_{k+1}+1$ is odd, thus
$t_{k+2}\dots t_n=0\dots 0$, and we retrieve the form prescribed by the first point of the proposition. 
\item  If $M$ is odd, then $M\neq b$ and $U_{k+1}+[t_{k+1}\neq0 {\rm\ and\ } t_{k+1} {\rm\ is\ even}]=U_{k+1}$ is even, thus
$t_{k+2}=\max\{s_1,\dots, s_k,t_{k+1}\}+1=M+1$, which is even. Next, we have 
$U_{k+1}+[t_{k+1}\neq0 {\rm\ and\ } t_{k+1} {\rm\ is\ even}]+ [t_{k+2}\neq0 {\rm\ and\ } t_{k+2} {\rm\ is\ even}]=U_{k+1}+1$ is odd, 
and this implies as above that $t_{k+3}\dots t_n=0\dots 0$, 
and we retrieve the second point of the proposition.
\end{itemize}

\noindent
For the case when $U_{k+1}$ is odd, in a similar way we have $t_{k+1}\dots t_n=0\dots 0$.

\qed

\medskip
The next proposition is the `first' counterpart of the previous one.

\begin{Pro}
\label{pro:pro_Rb_even2}
Let $b\geq2$ and even, $k\leq n-2$ and $\bsb s=s_1s_2\ldots s_k$. If $\bsb t$ is the $\preccdot$-first sequence in $\bsb s\,|\,R_n(b)$, then $\bsb t$ has one of the following forms:
\begin{enumerate}
\item $\bsb t=\bsb sM0\ldots0$ if $U_{k+1}$ is odd and  $M$ is even,
\item $\bsb t=\bsb sM(M+1)0\ldots0$ if $U_{k+1}$ is odd and $M$ is odd,
\item $\bsb t=\bsb s0\ldots0$ if $U_{k+1}$ is even,
\end{enumerate}
where $M=\min\{b,\max\{s_i\}_{i=1}^k+1\}$ and $U_{k+1}=\sum_{i=1}^k[s_i\neq0 {\rm\ and\ } s_i {\rm\ is\ even}]$.
%where $M=\max\{s_i\}_{i=1}^k+1$.
\end{Pro}

Based on Propositions \ref{pro:pro_Rb_even1}  and \ref{pro:pro_Rb_even2} we have the following theorem, its proof
is similar with that of Theorem \ref{k_odd}.
\begin{The} 
\label{k_even}
For any $n\geq 1$, $b\geq 2$ and even, $R_n(b)$ listed in $\preccdot$ order
is a $3$-adjacent Gray code.
\end{The}

It is worth to mention that, neither $\prec$ for even $b$, nor $\preccdot$ for odd $b$
yields a Gray code on $R_n(b)$.
Considering $b\geq n$ in Theorem \ref{k_odd} and \ref{k_even}, the bound $b$ does not actually provide any restriction, 
and in this case $R_n(b)=R_n$, and we have the following corollary.

\begin{Co} For any $n\geq 1$, $R_n$ listed in both $\prec$ and $\preccdot$ order
are $3$-adjacent Gray codes.
\end{Co}

\begin{The} 
\label{k}
For any $b\geq 1$ and odd, $n>b$, $R_n^*(b)$ listed in $\prec$ order
is a $5$-Gray code.
\end{The} 
\proof
For two integers $a$ and $b$, $0<a\leq b$, we define $\tau_{a,b}$ as the length $b-a$ increasing sequence $(a+1)(a+2) \ldots (b-1)b$, and 
$\tau_{a,b}$ is vanishingly empty if $a=b$. 
Imposing to a sequence $\bsb s$ in $R_n(b)$ to have its largest element equal to $b$ (so, to belong to $R^*_n(b)$)
implies that either $b$ occurs in $\bsb s$ before its last position, or $\bsb s$ ends with $b$, and in this case the 
tail of $\bsb s$ is $\tau_{a,b}$ for an appropriate $a<b$.
More precisely, in the latter case, $\bsb s$ has the form $s_1s_2\ldots s_j\tau_{a,b}$, for some $j$ and $a$, with 
$a=\max\{s_i\}_{i=1}^j$ and $j=n-(b-a)$.

Now let $\bsb s=s_1s_2\ldots s_n\prec\bsb t=t_1t_2\ldots t_n$ be two consecutive sequences in the $\prec$ ordered list for $R^*_n(b)$,
and let $k\leq n-3$
be the leftmost position where $\bsb s$ and $\bsb t$ differ, thus
$s_1s_2\ldots s_{k-1}=t_1t_2\ldots t_{k-1}$.
It follows that $\bsb s$ is the $\prec$-last sequence in $R^*_n(b)$ having the 
prefix $s_1s_2\ldots s_k$, and using Proposition \ref{pro:pro_Rb_odd1} and the notations therein, 
by imposing that $\max\{s_i\}_{i=1}^n$ is equal to $b$, we have:

\begin{itemize}
\item if $\sum_{i=1}^ks_i$ is odd, then $\bsb s$ has the form $s_1s_2\ldots s_k0\ldots0\,\tau_{a,b}$, where $a=\max\{s_i\}_{i=1}^k$,
\item if $\sum_{i=1}^ks_i$ is even,  then $\bsb s$ has one of the following forms:
   \begin{itemize}
   \item $s_1s_2\ldots s_k M0\ldots0\,\tau_{M,b}$, or
   \item $s_1s_2\ldots s_kM(M+1)0\ldots0\,\tau_{M+1,b}$.
   \end{itemize}
\end{itemize}
When the above $\tau$'s suffixes are empty, we retrieve precisely the three cases in Proposition \ref{pro:pro_Rb_odd1}.

Similarly, $\bsb t$ is the $\prec$-first sequence in $R^*_n(b)$ having the 
prefix $t_1t_2\ldots t_{k-1}t_k=s_1s_2\ldots s_{k-1}t_k$. 
Since by the definition of $\prec$ order we have that $t_k=s_k+1$ or $t_k=s_k-1$, it follows that
$\sum_{i=1}^kt_i$ and $\sum_{i=1}^ks_i$ have different parity (that is, $\sum_{i=1}^kt_i$ is odd if
and only if $\sum_{i=1}^ks_i$ is even), and by Proposition \ref{pro:pro_Rb_odd2} and replacing 
for notational convenience $M$ by $M'$, we have:

\begin{itemize}
\item if $\sum_{i=1}^ks_i$ is odd, then $\bsb t$ has the form $t_1t_2\ldots t_k0\ldots0\,\tau_{a',b}$, where $a'=\max\{t_i\}_{i=1}^k$,
\item if $\sum_{i=1}^ks_i$ is even,  then $\bsb t$ has one of the following forms:
   \begin{itemize}
   \item $t_1t_2\ldots t_k M'0\ldots0\,\tau_{M',b}$, or
   \item $t_1t_2\ldots t_kM'(M'+1)0\ldots0\,\tau_{M'+1,b}$.
   \end{itemize}
\end{itemize}

\noindent
With these notations, since $t_k\in \{s_k+1,s_k-1\}$, it follows that 
\begin{itemize}
\item if $\sum_{i=1}^ks_i$ is odd, then $a'\in\{a-1,a,a+1\}$, and so the length of $\tau_{a,b}$
and that of $\tau_{a',b}$ differ by at most one; and 
\item if $\sum_{i=1}^ks_i$ is even, then 
$M'\in\{M-1,M,M+1\}$, and the length of the non-zero tail of $\bsb s$ and that of $\bsb t$ 
(defined by means of $\tau$ sequences) differ by at most two.
\end{itemize}

Finally, the whole sequences $\bsb s$ and $\bsb t$ differ in at most five (not necessarily adjacent) positions,
and the statement holds.
\qed

\section{Generating algorithms}

An exhaustive generating algorithm is one generating 
all sequences in a combinatorial class, with some predefined properties 
({\it e.g.}, having the same length).
Such an algorithm is said to run in {\it constant amortized time} if it generates each object in $O(1)$ time, in amortized sense. 
In \cite{Ruskey} the author called such an algorithm {\it CAT algorithm} and shows that a recursive generating algorithm 
satisfying the following properties is precisely a CAT algorithm:
\begin{itemize}
\item Each recursive call either generates an object or produces at least two recursive calls;
\item The amount of computation in each recursive call is proportional to the degree of the call 
(that is, to the number of subsequent recursive calls produced by current call).
\end{itemize}  

Procedure {\sc Gen1} in Fig. \ref{fig:alg_Rnb_odd}
generates all sequences belonging to $R_n(b)$ in Reflected Gray Code Order. 
Especially when $b$ is odd, the generation induces a 3-adjacent Gray code. The bound $b$ and the generated sequence $s=s_1s_2\ldots s_n$ are global. The $k$ parameter is the position where the value is to be assigned (see line 8 and 13); the $dir$ parameter represents the direction of sequencing for $s_k$, whether it is up (when $dir$ is even, see line 7) or down (when $dir$ is odd, see line 12); and $m$ is such that $m+1$ is the 
the maximum  value that can be assigned to $s_k$, that is, 
$\min\{b-1,\max\{s_i\}_{i=1}^{k-1}\}$ (see line 5).

The algorithm initially sets $s_1=0$, and the recursive calls are triggered by the initial call {\sc Gen1}$(2,0,0)$. For the current position $k$, the algorithm assigns a value to $s_k$ (line 8 or 13) followed by recursive calls in line 10 or 15.
% A unique sequence of length $n$ is generated when $k=n+1$ (line 5).
This scheme guarantees that each recursive call will produce subsequent recursive calls until $k=n+1$ (line 4), 
that is, when a sequence of length $n$ is generated and printed out by {\sc Type()} procedure. This process eventually generates 
all sequences in $R_n(b)$. In addition, by construction, algorithm {\sc Gen1} satisfies the previous 
CAT desiderata, and so it is en efficient exhaustive generating algorithm.

\begin{figure}[h]
\centering
\begin{tabular}{|c|}
\hline
\begin{minipage}[c]{.46\linewidth}
\begin{tabbing}\hspace{0.8cm}\=\\
{\small 01} {\bf procedure} {\sc Gen1}($k$, $dir$, $m$: integer)\\
{\small 02} {\bf global} $s$, $n$, $b$: integer;\\
{\small 03} {\bf local} $i$, $u$: integer;\\
{\small 04}  {\bf if} $k=n+1$ {\bf then} {\sc Type()};\\
{\small 05}  {\bf else}  \={\bf if} $m=b$ {\bf then} $m:=b-1$; {\bf endif}\\
{\small 06}              \>{\bf  if}  $dir$ mod $2=0$\\
{\small 07}              \>{\bf then} \={\bf for} \=$i:=0$ {\bf to} $m+1$ {\bf do}\\
{\small 08}              \>           \>          \>$s_k:=i$;\\
{\small 09}              \>           \>          \>{\bf if} $m<s_k$ {\bf then} $u:=s_k$; {\bf else} $u:=m$; {\bf endif}\\
{\small 10}              \>           \>          \>{\sc Gen1$(k+1,i,u)$};\\   
{\small 11}              \>           \>{\bf endfor}\\
{\small 12}              \>{\bf else} \>{\bf for} $i:=m+1$ {\bf downto} $0$ {\bf do}\\
{\small 13}              \>           \>          \>$s_k:=i$;\\
{\small 14}              \>           \>          \>{\bf if} $m<s_k$ {\bf then} $u:=s_k$; {\bf else} $u:=m$; {\bf endif}\\
{\small 15}              \>           \>          \>{\sc Gen1}$(k+1,i+1,u)$;\\	
{\small 16}              \>           \>{\bf endfor}\\
{\small 17}              \>{\bf endif}\\
{\small 18}  {\bf endif}\\
{\small 19}  {\bf end procedure.}
\end{tabbing}
\end{minipage}
\\
\hline
\end{tabular}
\caption{
Reflected Gray Code Order generating algorithm for $R_n(b)$;
it produces a 3-Gray code when $b$ is odd.
\label{fig:alg_Rnb_odd}}
\end{figure}

Similarly, the call {\sc Gen2}$(2,0,0)$ of the algorithm in Fig. \ref{fig:alg_Rnb_even} generates sequences in $R_n(b)$ in 
co-Reflected Gray Code Order, and in particular when $b$ is even, a 3-adjacent Gray code for these sequences. 
And again it satisfies the CAT desiderata, and so it is en efficient exhaustive generating algorithm.

Finally, algorithm {\sc Gen3} in Fig. \ref{fig:alg_Pnb_odd} generates the set $R^*_n(b)$ in Reflected Gray Code Order
and produces a 5-Gray code if $b$ is odd. It mimes algorithm {\sc Gen1} and the only differences consist in an additional
parameter $a$ and lines 5, 6, 13 and 19, and its main call is {\sc Gen3}$(2,0,0,0)$. 
Parameter $a$ keeps track of the maximum value in the prefix $s_1s_2\dots s_{k-1}$
of the currently generated sequence, and it is updated in lines 13 and 19.
Furthermore, 
when the current position $k$ belongs to a $\tau$-tail (see the proof of Theorem \ref{k}),
that is, condition $k=n+1+a-b$ in line 5 is satisfied, then the imposed value is written in this position, 
and similarly for the next two positions. Theorem \ref{k} ensures that there are no differences between the current 
sequence and the previous generated one beyond position $k+2$, and thus a new sequence in $R^*_n(b)$ is generated.
And as previously, {\sc Gen3} is a CAT generating algorithm.

\begin{figure}[!h]
\centering
\begin{tabular}{|c|}
\hline

\begin{minipage}[c]{.46\linewidth}
\begin{tabbing}\hspace{0.8cm}\=\\
{\bf procedure} {\sc Gen2}($k$, $dir$, $m$: integer)\\
{\bf global} $s$, $n$, $b$: integer;\\
{\bf local} $i$, $u$: integer;\\
{\bf if} $k=n+1$ {\bf then} {\sc Type()};\\
{\bf else} \= {\bf if} $m=b$ {\bf then} $m:=b-1$;  {\bf endif}\\
           \> {\bf  if} $dir$ mod 2$=0$ \\
           \> {\bf then} \={\bf for} \= $i:=0$ {\bf to} $m+1$ {\bf do}\\
           \>           \>           \>$s_k:=i$;\\
           \>           \>           \>{\bf if} $m<s_k$ {\bf then} $u:=s_k$; {\bf else} $u:=m$; {\bf endif}\\
           \>           \>           \>{\bf if} $s_k=0$ \={\bf then} {\sc Gen2}$(k+1, 0,u)$;\\
           \>           \>           \>                 \>{\bf else} {\sc Gen2}$(k+1,  i+1,u)$;\\    
           \>           \>           \>{\bf endif}\\
           \>           \>{\bf endfor}\\
           \> {\bf else}\>{\bf for}\>$i:=m+1$ {\bf downto} $0$ {\bf do}\\ 
           \>           \>         \>$s_k:=i$;\\
           \>           \>         \>{\bf if} $m<s_k$  {\bf then} $u:=s_k$; {\bf else} $u:=m$; {\bf endif}\\
           \>           \>         \>{\bf if } $s_k=0$ \={\bf then} {\sc Gen2}$(k+1,  1,u)$;\\
           \>           \>         \>                  \>{\bf else} {\sc Gen2}$(k+1,  i,u)$;\\
           \>           \>         \>{\bf endif}\\
           \>           \>{\bf endfor}\\
           \> {\bf endif}\\
{\bf endif}\\
{\bf end procedure.}
\end{tabbing}
\end{minipage}
\\
\hline
\end{tabular}
\caption{
Co-Reflected Gray Code Order generating algorithm for $R_n(b)$;
it produces a 3-Gray code when $b$ is even.
\label{fig:alg_Rnb_even}}
\end{figure}

\begin{figure}[!h]
\centering
\begin{tabular}{|c|}
\hline
\begin{minipage}[c]{.46\linewidth}
\begin{tabbing}\hspace{0.8cm}\=\\
{\small 01} {\bf procedure} {\sc Gen3}($k$, $dir$, $m$, $a$: integer)\\
{\small 02} {\bf global} $s$, $n$, $b$: integer;\\
{\small 03} {\bf local} $i$, $u$, $\ell$: integer;\\
{\small 04} {\bf if} $k=n+1$ {\bf then} {\sc Type()};\\
{\small 05} {\bf else} \= {\bf if} $k=n+1 + a-b$\\
{\small 06}            \> {\bf then} \= {\bf for} $i:=0$ {\bf to} $2$ {\bf do} {\bf if} $k+i\leq n$ {\bf then} $s_{k+i}:=a+1+i$; {\bf endif} {\bf endfor}\\
{\small 07}            \>            \> {\sc Type()};\\
{\small 08}            \> {\bf else} \={\bf if} $m=b$ {\bf then} $m:=b-1$; {\bf endif}\\
{\small 09}            \>           \>{\bf  if} $dir$ mod 2$=0$\\
{\small 10}            \>           \>{\bf then} \={\bf for} \= $i:=0$ {\bf to} $m+1$ {\bf do}\\
{\small 11}            \>           \>           \>          \>$s_k:=i$;\\
{\small 12}            \>           \>           \>          \>{\bf if} $m<s_k$   {\bf then} $u:=s_k$; {\bf else} $u:=m$; {\bf endif}\\
{\small 13}            \>           \>           \>          \>{\bf if} $a<s_k$ {\bf then} $\ell:=s_k$; {\bf else} $\ell:=a$; {\bf endif}\\
{\small 14}            \>           \>           \>          \>{\sc Gen3$(k+1,i,u,\ell)$};\\    
{\small 15}            \>           \> \>{\bf endfor}\\         
{\small 16}            \>           \>{\bf else} \>{\bf for} \>$i:=m+1$ {\bf downto} $0$ {\bf do}\\
{\small 17}            \>           \>           \>          \>$s_k:=i$;\\
{\small 18}            \>           \>           \>          \>{\bf if} $m<s_k$ {\bf then} $u:=s_k$; {\bf else} $u:=m$; {\bf endif}\\
{\small 19}            \>          \>           \>          \>{\bf if} $a<s_k$ {\bf then} $\ell:=s_k$; {\bf else} $\ell:=a$; {\bf endif}\\
{\small 20}            \>           \>           \>          \>{\sc Gen3}$(k+1,i+1,u,\ell)$;\\	
{\small 21}            \>           \> \>{\bf endfor}\\   
{\small 22}            \>           \>{\bf endif}\\
{\small 23}            \>{\bf endif}\\
{\small 24} {\bf endif}\\
{\small 25} {\bf end procedure.}
\end{tabbing}	
\end{minipage}
\\
\hline
\end{tabular}
\caption{
Generating algorithm for $R_n^*(b)$, $n>b\geq 1$, with respect to
Reflected Gray Code Order; it produces a 5-Gray code when $b$ is odd.
\label{fig:alg_Pnb_odd}}
\end{figure}

\medskip

\noindent
{\bf Final remarks.}
We suspect that the upper bounds
3 in Theorems \ref{k_odd} and \ref{k_even}, and 5 in Theorem \ref{k} are not tight, and a natural question
arises: are there more restrictive Gray codes for $R_n(b)$ and for $R^*_n(b)$ with $b$ odd?
Finally, is there a natural order relation inducing a Gray code on $R^*_n(b)$ when $b$ is even?

\end{document}